\documentclass[12pt,reqno]{elsarticle}

\usepackage{amssymb}
\usepackage{amsthm}

\usepackage[colorlinks=true,linkcolor=black, citecolor=blue, urlcolor=blue]{hyperref}
\usepackage{pgffor}
\usepackage{epsfig}
\usepackage{graphics}
\usepackage{amsfonts}
\usepackage{amsmath}
\usepackage{latexsym,mathrsfs}
\usepackage[nottoc, notlof, notlot]{tocbibind}
\newtheorem{theorem}{Theorem}[section]
\newtheorem{lemma}{Lemma}[section]
\newtheorem{proposition}{Proposition}[section]
\newtheorem{corollary}{Corollary}[section]

\newtheorem{definition}{Definition}[section]

\newtheorem{remarks}{Remarks}[section]

\newtheorem{remark}{Remark}[section]

\newtheorem*{hyp}{}

\DeclareMathOperator{\supp}{supp} 
 \numberwithin{equation}{section}

\journal{}
\begin{document}

\begin{frontmatter}


\title{\textbf{On the projections of mutual multifractal spectra} }

\author{DOUZI Zied and SELMI Bilel \footnote{Faculty of sciences of Monastir,
Department of mathematics, 5000-Monastir, Tunisia.\\
{E-mail addresse:\;zied.douzi@fsm.rnu.tn,\;\;bilel.selmi@fsm.rnu.tn
}}}



\date{May 20, 2011}
\begin{abstract}
The aim of this article is  to study the behaviour of the relative
multifractal spectrum under projections. First of all, we depict a
relationship between the mutual multifractal spectra of a couple of
measures $(\mu,\nu)$ and its orthogonal projections in Euclidean
space. As an application, we improve Svetova's result (Tr.
Petrozavodsk. Gos. Univ. Ser. Mat., {\bf 11} (2004), 41-46) and
study the mutual multifractal analysis of the projections of
measures.

\end{abstract}
\end{frontmatter}

\addcontentsline{toc}{section}{On the projections of mutual
multifractal spectra}
\section{Introduction}

In the previous years, there has been great interest in
understanding the fractal dimensions of projections of sets and
measures. The first significant work in this area was the result of
Marstrand \cite{21} who showed a well-known theorem according to
which the Hausdorff dimension of a planar set is preserved under
orthogonal projections. In \cite{19}, Kaufman had employed potential
theoretic methods in order to prove Marstrand result, which has been
generalized later by Mattila in \cite{22}. Let us mention that
Falconer et al \cite{12,13} have proved that the packing dimension
of the projected set or measure will be the same for almost all
projections. Other works were carried out in this sense for classes
of similar measures in euclidean and symbolic spaces \cite{BF, 15,
16, 31, 32}. However, despite these substantial advances for fractal
sets, only very little is known about the multifractal structure of
projections of measures \cite{1, 4, 28, B, SB, SBS}.

Based on some ideas of multifractal formalism given by Olsen and
Peyri\`{e}re \cite{25, JP}, Svetova introduced in \cite{5, 6, 7, 8}
a new formalism for a multifractal analysis of one measure with
respect to an other. This formalism is called by the mutual
multifractal formalism and for which Svetova studied some basic
properties. More specifically, given two compactly supported Borel
probability measures $\mu$ and $\nu$ on $\mathbb{R}^n$ and
$\alpha,\beta\in\mathbb{R}$, Svetova estimated the size of the
iso-H\"{o}lder set
$$
E_{\mu,\nu}(\alpha,\beta) =  \Big\{ x \in \supp \mu \cap \supp \nu ;
\;\;{\alpha}_{\mu}(x) =\alpha
\quad\text{and}\quad{\alpha}_{\nu}(x)=\beta\Big\},
$$
where ${\alpha}_{\mu}(x)=\lim_{r\to 0}\frac{\log \mu(B(x,r))}{\log
r}$ and $B(x,r)$ is the closed ball of center $x$ and radius $r$.
The mutual multifractal analysis of a measures allows to relate the
Hausdorff and packing dimensions of these levels sets to the
Legendre transform of some multifractal functions. There has
recently been a great interest for this subject and positive results
have been written in various situations in the dynamic contexts
\cite{BSB, BP, BP1, BP2, BCW, Ol2}. Recently, many authors were
interested in mutual (mixed) multifractal spectra, see for example
\cite{D, D11, D21, D31, D41, D51, LOW, Ol1, Ol, SSSS}. We write for
$\gamma\geq0$,
 $$
\mathcal{B}_{\mu,\nu}(\gamma)=\left\{x\in \supp\mu\cap \supp\nu ;\:
\:\lim_{r\to 0}\frac{\log \big(\mu(B(x,r))\big)}{\log\big(
\nu(B(x,r)\big)}=\gamma\right\}.
 $$
It is clear that
 $$
\displaystyle
\bigcup_{\displaystyle{\substack{(\alpha,\beta)\in\mathbb{R}_+\times
\mathbb{R}_+^*, \\\frac{\alpha}{\beta}=\gamma}}}\;
E_{\mu,\nu}\big(\alpha,\beta\big) \subset
\mathcal{B}_{\mu,\nu}(\gamma).
 $$
The latter union is composed by an uncountable number of pairwise
disjoint nonempty sets. Then, the Hausdorff and packing dimension of
$\mathcal{B}_{\mu,\nu}(\gamma)$ is fully carried by some subset
$E_{\mu,\nu}\big(\alpha,\beta\big)$.  Also, Selmi et al.
investigated the projection properties of the $\nu$-Hausdorff, and
the $\nu$-packing dimensions of $\mathcal{B}_{\mu,\nu}(\gamma)$  in
\cite{4}. In this article, they derived global bounds on the
relative multifractal dimensions of a projection of a measures in
terms of its original relative multifractal dimensions. It is more
difficult to obtain a lower and upper bound for the dimension of the
set $\mathcal{B}_{\mu_V,\nu_V}(\gamma)$, where $V$ is a linear
subspace of $\mathbb{R}^n$.

The purpose of this paper is to improve Svetova's result and to
propose a sufficient condition that gives the lower bound for the
Hausdorff and the packing dimensions of
$\mathcal{B}_{\mu_V,\nu_V}(\gamma)$. Our first aim is to study the
behavior of the multual Hausdorff, packing and pre-packing
dimensions under projections. The second aim is to investigate a
relationship between the mutual multifractal spectra and its
projections onto a lower dimensional linear subspace.

\section{Preliminaries and Results}
Let us recall the multifractal formalism introduced by Svetova in
\cite{6}. Let $\mu$ and $\nu$ be two compactly supported Borel
probability measures on $\mathbb{R}^n$. We denote by $\supp\mu$ the
topological support of $\mu$.

\begin{definition}
For $q, t, s\in\mathbb{R}$, $E \subseteq{\mathbb R}^n$ and
$\delta>0$, we define

$$\overline{{\mathcal P}}^{q,t,s}_{\mu,\nu,\delta}(E) =\displaystyle
 \sup \sum_i \mu(B(x_i,r_i))^q\nu (B(x_i, r_i))^t(2r_i)^s,$$
where the supremum is taken over all  centered $\delta$-packing of
$E$, $$\overline{{\mathcal P}}^{q,t,s}_{\mu,\nu}(E)
=\displaystyle\inf_{\delta>0}\overline{{\mathcal
P}}^{q,t,s}_{\mu,\nu,\delta}(E),$$ and we introduce the  generalized
packing measure relatively to $\mu$ and $\nu$
 $$
{\mathcal P}^{q,t,s}_{\mu,\nu}(E) = \inf_{E \subseteq
\bigcup_{i}E_i} \sum_i \overline{\mathcal P}^{q,t,s}_{\mu,\nu}(E_i).
 $$
In a similar way we define the  generalized Hausdorff measure
relatively to $\mu$ and $\nu$ by \\For $E \subseteq{\mathbb R}^n$
and $\delta>0$,
 $$
\overline{{\mathcal H}}^{q,t,s}_{\mu,\nu,\delta}(E) =
\displaystyle\inf \sum_i \mu(B(x_i,r_i))^q\nu (B(x_i,
r_i))^t(2r_i)^s,
 $$
where the infinimum is taken over all centered $\delta$-covering  of
$E$,
 $$
 \overline{{\mathcal H}}^{q,t,s}_{\mu,\nu}(E) =
\displaystyle\sup_{\delta>0}\overline{{\mathcal
H}}^{q,t,s}_{\mu,\nu,\delta}(E),
 $$
and we introduce the  generalized Hausdorff measure relatively to
$\mu$ and $\nu$
 $$
{\mathcal H}^{q,t,s}_{\mu,\nu}(E)=\displaystyle\sup_{F\subseteq
E}\overline{{\mathcal H}}^{q,t,s}_{\mu,\nu}(F),
 $$
with the conventions $0^q = \infty$ for $q\leq0$ and $0^q = 0$ for
$q>0$.
\end{definition}
\begin{remarks}\par\noindent
\begin{enumerate}
\item The functions ${\mathcal H}^{q,t,s}_{\mu,\nu}$ and ${\mathcal
P}^{q,t,s}_{\mu,\nu}$ are metric outer measures and thus measures on
the Borel family of subsets of $\mathbb{R}^n$. An important feature
of the Hausdorff and packing measures is that ${\mathcal
P}^{q,t,s}_{\mu,\nu}\leq{\overline{\mathcal P}}^{q,t,s}_{\mu,\nu}$
and that there exists an integer $\xi\in\mathbb{N}$, such that
${\mathcal H}^{q,t,s}_{\mu,\nu}\leq\xi{\mathcal
P}^{q,t,s}_{\mu,\nu}$ (see \cite{S}).

\item In the special case where $q=0$ or $t=0$, the mutual multifractal spectra is strictly
related to Olsen's multifractal formalism \cite{25}.

\item The mutual multifractal spectra represents the relative multifractal
analysis introduced by Cole \cite{10} in the case where $s=0$. Other
works were carried out in this sense in probability and symbolic
spaces \cite{NBC, D1, D3, D4}.
\end{enumerate}
\end{remarks}
\begin{proposition}(\cite{6, S})
\par\noindent\begin{enumerate}\item  There exists a unique number
$b_{\mu,\nu}^{q,t}(E)\in[-\infty,+\infty]$ such that
 $$
{\mathcal H}^{q,t,s}_{\mu,\nu}(E)=\left\{\begin{matrix}
 \infty &\text{if}& s < b_{\mu,\nu}^{q,t}(E),\\ \\
 0 & \text{if}&  b_{\mu,\nu}^{q,t}(E) < s.
 \end{matrix}\right.
 $$
\item  There exists a unique number
$B_{\mu,\nu}^{q,t}(E)\in[-\infty,+\infty]$ such that
 $$
{\mathcal P}^{q,t,s}_{\mu,\nu}(E)=\left\{\begin{matrix}
 \infty &\text{if}& s < B_{\mu,\nu}^{q,t}(E),\\
 \\
 0 & \text{if}&  B_{\mu,\nu}^{q,t}(E) < s.\end{matrix}\right.
 $$
\item  There exists a unique number
$\Lambda_{\mu,\nu}^{q,t}(E)\in[-\infty,+\infty]$ such that
 $$
{\overline{\mathcal P}}^{q,t,s}_{\mu,\nu}(E)=\left\{\begin{matrix}
\infty &\text{if}& s < \Lambda_{\mu,\nu}^{q,t}(E),\\
\\
0 & \text{if}& \Lambda_{\mu,\nu}^{q,t}(E) < s.\end{matrix}\right.
 $$
\end{enumerate}
\end{proposition}

Let $E\subseteq \mathbb{R}^n$ and $q,t\in\mathbb{R}.$ We can remark
that \begin{center}$b^{q,t}_{\mu,\nu}(E)\leq B^{q,t}_{\mu,\nu}(E)
\leq\Lambda^{q,t}_{\mu,\nu}(E).$
\end{center}
Then we are able to define the multifractal dimension functions
$b_{\mu,\nu},$ $B_{\mu,\nu}$ and $\Lambda_{\mu,\nu}$:
$\mathbb{R}^2\rightarrow [-\infty, +\infty]$ by
 \begin{center}$
b_{\mu,\nu}(q,t)=b_{\mu,\nu}^{q,t}(\supp\mu\cap \supp\nu), \quad
B_{\mu,\nu}(q,t)=B_{\mu,\nu}^{q,t}(\supp\mu\cap \supp\nu)$
\end{center}
\begin{center}\;\text{and}\;
$\Lambda_{\mu,\nu}(q,t)=\Lambda_{\mu,\nu}^{q,t}(\supp\mu\cap
\supp\nu).
 $\end{center}
It is well known that the functions $b_{\mu,\nu}$, $B_{\mu,\nu}$ and
$\Lambda_{\mu,\nu}$ are decreasing and $B_{\mu,\nu}$,
$\Lambda_{\mu,\nu}$ are convex (see \cite{S}).

\section{Projection results}
Let $m$ be an integer with $0<m<n$ and $G_{n,m}$ stand for the
Grassmannian  manifold of all $m$-dimensional linear subspaces of
$\mathbb{R}^n$. We denote by $\gamma_{n,m}$ the invariant Haar
measure on $G_{n,m}$ such that $\gamma_{n,m}(G_{n,m})=1$. For $V\in
G_{n,m}$, we define the projection map, $\pi_V:
\mathbb{R}^n\longrightarrow V$ as the usual orthogonal projection
onto $V$. Now, for a Borel probability measure $\mu$ on
$\mathbb{R}^n$, supported on the compact set $\supp\mu$ and for
$V\in G_{n,m}$, we define $\mu_V$, the projection of $\mu$ onto $V$,
by
 $$
\mu_V(A)=\mu(\pi_V^{-1}(A)),\quad \forall A\subseteq V.
 $$
Since $\mu$ has a compact support, $\supp\mu_V=\pi_V(\supp\mu)$ for
all $V\in G_{n,m}$.

 \bigskip
In the following, we are interested about the behavior of mutual
Hausdorff, packing and pre-packing dimensions under projections.
Throughout this paper, we suppose that $\supp\mu=\supp\nu$. We are
based on ideas of Selmi et al in \cite{4}, to show the following
results.

\begin{theorem}\label{g}
Let $\mu$ and $\nu$ be two compactly supported Borel probability
measures on $\mathbb{R}^n$ and $E\subseteq \supp\mu$. Then, for
$(q,t)\in\big(]-\infty,0]^2\big)\cup\big(]-\infty,0]\times
[0,1]\big)\cup\big([0,1]\times ]-\infty,0]\big)$ and for all $V\in
G_{n,m},$ we have
 $$
\Lambda_{\mu_V,\nu_V}^{q,t}(E)\leq \Lambda_{\mu,\nu}^{q,t}(E).
 $$
\end{theorem}

\noindent{\bf Proof.} Let $s\in\mathbb{R}$ such that
$\Lambda_{\mu,\nu}^{q,t}(E)<s.$ Consider $V\in G_{n,m}$ and fix
$\delta>0.$
\\Let $\Big(B_i=B(x_i,r_i)\Big)_i$ be a $\delta$-centered packing of $\pi_V(E)$.
There exists an integer $K_m$ depending on $m$ only such that we can
divide up the balls $B(x_{i},2r_i)$ into $K\le K_m$ families of
disjoint balls $\mathcal{B}_1,\dots,\mathcal{B}_K$. Let $1\le l\le
K$. For each $B(x_i,r_i)\in \mathcal{B}_l$, denote $E_i=E\cap
\pi_V^{-1}\big(B(x_i,r_i)\big)$. \\We have $E_i\subset\bigcup_{y\in
E_i}B(y,r_i)$, so Besicovitch's covering theorem \cite{23} provides
a positive integer $K_n$ as well as $K_i\le K_n$ families of
pairwise disjoint balls ${\mathcal
B_{i,k}}=\left\{B^{(i,k)}_j=B(y^{(i,k)}_j,r_{ijk});  \:
r_{ijk}=\displaystyle\frac {r_i} 2\right\}$, $1\le k\le K_i$,
extracted from $\Big\{B(y,r_i)\Big\}_{y\in  E_i}$ such that
 $$
E_i\subset \bigcup_{k=1}^{K_i}\bigcup_{j}B^{(i,k)}_j.
 $$
$\bullet\:\underline{Case}\; 1:$ For $q\leq 0$ and $t\le0$, we have
\begin{eqnarray*}
\sum_i\mu_V(B_i)^q\nu_V(B_i)^t(2r_i)^s&\leq&
2^s\sum_i\mu\big(B_j^{(i,k)}\big)^q\nu \big
(B_j^{(i,k)}\big )^t(2r_{ijk})^s\\
&\leq&2^s\sum_{i,j}\sum_{k=1}^{k_i}\mu\big(B_j^{(i,k)}\big)^q
\nu\big(B_j^{(i,k)}\big)^t(2r_{ijk})^s.
\end{eqnarray*}
$\bullet\:\underline{Case}\; 2:$ For $q\leq 0$ and $0\leq t\leq 1$,
we have
\begin{eqnarray*}
\sum_i\mu_V(B_i)^q\nu_V(B_i)^t(2r_i)^s&\leq&
2^s\sum_i\mu\big(B_j^{(i,k)}\big)^q\nu \left
(\bigcup_{k=1}^{K_i}\bigcup_j B_j^{(i,k)} \right )^t(2r_{ijk})^s\\
&\leq&2^s\sum_{i,j}\sum_{k=1}^{k_i}\mu\big(B_j^{(i,k)}\big)^q
\nu\big(B_j^{(i,k)}\big)^t(2r_{ijk})^s.
\end{eqnarray*}
$\bullet\:\underline{Case}\;3:$ For $0\leq q \leq 1$ and $t\leq 0$,
we have
\begin{eqnarray*}
\sum_i\mu_V(B_i)^q\nu_V(B_i)^t(2r_i)^s&\leq&
2^s\sum_i\mu\left(\bigcup_{k=1}^{K_i}\bigcup_j
B_j^{(i,k)}\right)^q\nu \big
(B_j^{(i,k)} \big )^t(2r_{ijk})^s\\
&\leq&2^s\sum_{i,j}\sum_{k=1}^{k_i}\mu\big(B_j^{(i,k)}\big)^q
\nu\big(B_j^{(i,k)}\big)^t(2r_{ijk})^s.
\end{eqnarray*}

In all cases and  by construction,  since the balls $B(x_i,2r_i)\in
\mathcal{B}_l$ are pairwise disjoint, if $B(y,r)\in {\mathcal
B_{i,k}}$ and $B(y',r')\in {\mathcal B_{i',k'}}$ with $i\neq i'$,
then $B(y,r)\cap B(y',r')=\emptyset$. Consequently, we can collect
the balls $B(y,r)$ invoked in the above sum into at most $K_n$
centered packing of $E$. This holds for all $1\le l\le K$, so
 $$
\sum_i \mu_V(B_i)^q\nu_V(B_i)^t(2r_i)^s\le 2^sK_mK_n \sup\left
\{\sum_j \mu(B(y_j,r_j))^q\nu(B(y_j,r_j))^t(2r_j)^s\right\},
 $$
where the supremum is taken over all centered packing of $E$ by
closed balls of radius  $r$. Thus
 $$
\overline{\mathcal{P}}_{\mu_V,\nu_V,\delta}^{q,t,s}(\pi_V(E))\leq
2^sK_nK_m\overline{\mathcal{P}}_{\mu,\nu,\delta}^{q,t,s}(E).
 $$
Letting $\delta\downarrow0$, we obtain
\begin{equation}\label{h}
\overline{\mathcal{P}}_{\mu_V,\nu_V}^{q,t,s}(\pi_V(E))\leq2^s
K_nK_m\overline{\mathcal{P}}_{\mu,\nu}^{q,t,s}(E),
\end{equation}
and the result yields.
\begin{corollary}
Let $\mu$ and $\nu$ be two compactly supported Borel probability
measures on $\mathbb{R}^n$. Then for
$(q,t)\in\big(]-\infty,0]^2\big)\cup\big(]-\infty,0]\times
[0,1]\big)\cup\big([0,1]\times ]-\infty,0]\big)$ and for all $V\in
G_{n,m},$ we have
 $$
\Lambda_{\mu_V,\nu_V}(q,t)\leq \Lambda_{\mu,\nu}(q,t).
 $$ \end{corollary}

\noindent{\bf Proof.} The proof is a consequence of
\autoref{g}.$\hfill\square$

\begin{theorem}\label{j}
Let $\mu$ and $\nu$ be two compactly supported Borel probability
measures on $\mathbb{R}^n$. Then for
$(q,t)\in\big(]-\infty,0]^2\big)\cup\big(]-\infty,0]\times
[0,1]\big)\cup\big([0,1]\times ]-\infty,0]\big)$ and for all $V\in
G_{n,m},$ we have
 $$
B_{\mu_V,\nu_V}(q,t)\leq B_{\mu,\nu}(q,t).
 $$
\end{theorem}

\noindent{\bf Proof.} Let $s\in\mathbb{R}$ such that
$B_{\mu,\nu}(q,t)< s.$
\\Consider $F\subseteq\mathbb{R}^n$ and $V\in G_{n,m}$. Due to inequality \eqref{h}, we have
 $$
\overline{\mathcal{P}}_{\mu_V,\nu_V}^{q,t,s}(\pi_V(F))\leq
2^sK_nK_m\overline{\mathcal{P}}_{\mu,\nu}^{q,t,s}(F).
 $$
Since $\mathcal{P}_{\mu,\nu}^{q,t,s}(\supp\mu)=0$, there exists
$(E_i)_i$ a covering of $\supp\mu$ such that
$$\displaystyle\sum_i\overline{\mathcal{P}}_{\mu,\nu}^{q,t,s}(E_i)<1.$$
So, $\pi_V(\supp\mu)\subseteq\displaystyle\bigcup_i\pi_V(E_i)$ and we
have
\begin{eqnarray*}
\mathcal{P}_{\mu_V,\nu_V}^{q,t,s}(\supp\mu_V) &\leq&
\displaystyle\sum_i
\overline{\mathcal{P}}_{\mu_V,\nu_V}^{q,t,s}(\pi_V(E_i)) \\
&\leq&2^sK_nK_m\sum_i\overline{\mathcal{P}}_{\mu,\nu}^{q,t,s}(E_i)<\infty.
\end{eqnarray*}
Thus $B_{\mu_V,\nu_V}(q,t)\leq s$. $\hfill\square$
\begin{theorem} \label{k} Let $\mu$, $\nu$ be two compactly
supported Borel probability measures on $\mathbb{R}^n$. Then for
$(q,t)\in\big(]-\infty,0[^2\big)\cup\big(]-\infty,0[\times
]0,1]\big)\cup\big(]0,1]\times ]-\infty,0[\big)$ and for all $V\in
G_{n,m},$ we have
 $$
b_{\mu_V,\nu_V}(q,t)=b_{\mu,\nu}(q,t).
 $$
\end{theorem}
\noindent{\bf Proof.} Let's prove that $ b_{\mu,\nu}(q,t)\leq
b_{\mu_V,\nu_V}(q,t).$\\ Let $s\in\mathbb{R}$ such that
$s<b_{\mu,\nu}(q,t).$ Choose $F\subseteq\supp\mu$ and $V\in
G_{n,m}$. Fix $\delta>0$ and let $\Big(B_i=B(x_i,r_i)\Big)_i$ be a
$\delta$-centered covering of $F$. Let $E_i$ such that
$\pi_V^{-1}(E_i)=F\cap B(x_i,r_i)$. We have $E_i\subset\bigcup_{y\in
E_i}B(y,r_i)$, so Besicovitch's covering theorem provides a positive
integer $K_n$ as well as $K_i\le K_n$ families of pairwise disjoint
balls ${\mathcal B_{i,k}}=\Big\{B^{(i,k)}_j=B(y^{(i,k)}_j,r_{ijk});
\:r_{ijk}=\frac {r_i} 2\Big\}$, $1\le k\le K_i$, extracted from
$\Big\{B(y,r_i)\Big\}_{y\in  E_i}$ and such that
 $$
E_i\subset \bigcup_{k=1}^{K_i}\bigcup_{j}B^{(i,k)}_j.
 $$
$\bullet\:\underline{Case}\;1:$ For $q< 0$ and $t<0$, we have
\begin{eqnarray*}
\sum_i\mu(B_i)^q\nu(B_i)^t(2r_i)^s&\leq&2^s\sum_i\mu_V\big(B_j^{(i,k)}\big)^q
\nu_V\big(B_j^{(i,k)}\big)^t(2r_{ijk})^s\\
&\leq& 2^s\sum_{i,j}\sum_{k=1}^{k_i}\mu_V(B_j^{(i,k)})^q
\nu_V(B_j^{(i,k)})^t(2r_{ijk})^s.
\end{eqnarray*}
$\bullet\:\underline{Case}\;2:$ For $q< 0$ and $0< t\le1$, we have
\begin{eqnarray*}
\sum_i\mu(B_i)^q\nu(B_i)^t(2r_i)^s&\leq&2^s
\sum_i\mu_V(B_j^{(i,k)})^q\nu_V \left(\displaystyle \bigcup
_{k=1}^{k_i}\bigcup_jB_j^{(i,k)}\right)^t(2r_{ijk})^s\\&\leq&
2^s\sum_{i,j}\sum_{k=1}^{k_i}\mu_V(B_j^{(i,k)})^q
\nu_V(B_j^{(i,k)})^t(2r_{ijk})^s.
\end{eqnarray*}
$\bullet\:\underline{Case}\;3:$ For $0< q\le 1$ and $t<0$, we have
\begin{eqnarray*}
\sum_i\mu(B_i)^q\nu(B_i)^t(2r_i)^s&\leq&2^s \sum_i\mu_V
\left(\displaystyle\bigcup _{k=1}^{k_i}\bigcup_jB_j^{(i,k)}\right)^q
\nu_V(B_j^{(i,k)})^t(2r_{ijk})^s\\&\leq&
2^s\sum_{i,j}\sum_{k=1}^{k_i}\mu_V(B_j^{(i,k)})^q
\nu_V(B_j^{(i,k)})^t(2r_{ijk})^s.
\end{eqnarray*}
Thus
 $$
\overline{\mathcal{H}}_{\mu,\nu,\delta}^{q,t,s}(F)\leq 2^s
\overline{\mathcal{H}}_{\mu_V,\nu_V,\delta}^{q,t,s}(\pi_V(F)).
 $$
Letting $\delta\downarrow0$, we obtain
 $$
\overline{\mathcal{H}}_{\mu,\nu}^{q,t,s}(F)\leq2^s
\overline{\mathcal{H}}_{\mu_V,\nu_V}^{q,t,s}(\pi_V(F)).
 $$
Thus
\begin{eqnarray*}
\overline{\mathcal{H}}_{\mu,\nu}^{q,t,s}(F) &\leq&2^s
\overline{\mathcal{H}}_{\mu_V,\nu_V}^{q,t,s}(\pi_V(F)) \\
&\leq&2^s\overline{\mathcal{H}}_{\mu_V,\nu_V}^{q,t,s}(\pi_V(\supp\mu))\\
&\leq&2^s\mathcal{H}_{\mu_V,\nu_V}^{q,t,s}(\supp\mu_V).
\end{eqnarray*}
The arbitrary on $F$ implies that
\begin{equation} \label{www}
\mathcal{H}_{\mu,\nu}^{q,t,s}(\supp\mu)\leq2^s
\mathcal{H}_{\mu_V,\nu_V}^{q,t,s}(\supp\mu_V)
\end{equation}
and the result holds.

\bigskip
In order to prove the other inequality, let $E\subseteq\mathbb{R}^n$
and $s\in\mathbb{R}$ such that $b_{\mu,\nu}^{q,t}(E)<s.$ Fix $V\in
G_{n,m}$, $\delta>0$ and suppose that $\Big(B(x_i,r_i)\Big)_i$ is a
$\delta$-cover of $\pi_V(E)$. Denote
$E_i=E\bigcap\pi_V^{-1}\big(B(x_i,r_i)\big)$. We have
$E_i=\bigcup_{y\in E_i\displaystyle\cap\pi_V^{-1}(\{x_i\})}B(y,\frac
{r_i} n)$. By applying Besicovitch covering theorem, we find an
integer $K_n$, depending only on n as well as $K_i\leq K_n$ families
of pairwise disjoint balls
$\mathcal{B}_{i,k}=\Big\{B^{(i,k)}_j=B(y_j^{(i,k)},\frac {r_{ijk}}
n);\; r_{ijk}=\frac {r_i} 2 \Big\}$, $1\leq k\leq K_i$ such that
$$E\cap\pi_V^{-1}\big(B(x_i,r_i)\big)\subseteq\displaystyle
\bigcup_{k=1}^{K_i}\bigcup_{j}B^{(i,k)}_j.$$
$\bullet\;\underline{Case}\;1:$ For $q< 0$ and $t<0$, we have
\begin{eqnarray*}
\sum_i\mu_V(B_i)^q\nu_V(B_i)^t(2r_i)^s&\leq&2^s\sum_i\mu\big(B_j^{(i,k)}\big)^q
\nu\big(B_j^{(i,k)}\big)^t(2r_{ijk})^s\\&\leq&
2^s\sum_{i,j}\sum_{k=1}^{k_i}\mu(B_j^{(i,k)})^q
\nu(B_j^{(i,k)})^t(2r_{ijk})^s.
\end{eqnarray*}
$\bullet\:\underline{Case}\;2:$ For $q< 0$ and $0< t\le1$, we have
\begin{eqnarray*}
\sum_i\mu_V(B_i)^q\nu_V(B_i)^t(2r_i)^s&\leq&2^s
\sum_i\mu(B_j^{(i,k)})^q\nu\left(\displaystyle \bigcup
_{k=1}^{k_i}\bigcup_jB_j^{(i,k)}\right)^t(2r_{ijk})^s\\&\leq&
2^s\sum_{i,j}\sum_{k=1}^{k_i}\mu(B_j^{(i,k)})^q
\nu(B_j^{(i,k)})^t(2r_{ijk})^s.
\end{eqnarray*}
$\bullet\:\underline{Case}\;3:$ For $0< q\le 1$ and $t<0$, we have
\begin{eqnarray*}
\sum_i\mu_V(B_i)^q\nu_V(B_i)^t(2r_i)^s&\leq&2^s
\sum_i\mu\left(\displaystyle\bigcup
_{k=1}^{k_i}\bigcup_jB_j^{(i,k)}\right)^q
\nu(B_j^{(i,k)})^t(2r_{ijk})^s\\&\leq&
2^s\sum_{i,j}\sum_{k=1}^{k_i}\mu(B_j^{(i,k)})^q
\nu(B_j^{(i,k)})^t(2r_{ijk})^s.
\end{eqnarray*}

Then
 $$
\overline{\mathcal{H}}_{\mu_V,\nu_V,\delta}^{q,t,s}(\pi_V(E))\leq2^s
\overline{\mathcal{H}}_{\mu,\nu,\delta}^{q,t,s}(E).
 $$
Letting $\delta\downarrow0$, we obtain
$$
\overline{\mathcal{H}}_{\mu_V,\nu_V}^{q,t,s}(\pi_V(E))\leq2^s
\overline{\mathcal{H}}_{\mu,\nu}^{q,t,s}(E).
 $$
  Thus, given a subset $E$ of $\supp\mu$, $\pi_V(E)\subseteq \supp\mu_V$ and
\begin{eqnarray*}
   \overline{\mathcal{H}}_{\mu_V,\nu_V}^{q,t,s}(\pi_V(E)) &\leq&
   2^s\overline{\mathcal{H}}_{\mu,\nu}^{q,t,s}(E) \\
   &\leq& 2^s\mathcal{H}_{\mu,\nu}^{q,t,s}(\supp\mu).
\end{eqnarray*}
The arbitrary on $E$ implies that
\begin{center}$\mathcal{H}_{\mu_V,\nu_V}^{q,t,s}(\supp\mu_V)\leq
2^s\mathcal{H}_{\mu,\nu}^{q,t,s}(\supp\mu)$\end{center} and the
result holds. This achieves the proof of \autoref{k}.$\hfill\square$

\begin{theorem} Let $\mu$ and $\nu$ be two compactly
supported Borel probability measures on $\mathbb{R}^n$. Then for
$q,t\ge 1$ and all $V\in G_{n,m},$ we have
 $$
b_{\mu_V,\nu_V}({q,t})\ge b_{\mu,\nu}({q,t}).
 $$
\end{theorem}

\noindent{\bf Proof.} Fix $V\in G_{n,m}$ and $\delta>0$ and suppose
that $\big(B_i=B(x_i,r_i)\big)_i$ is a  $\delta$-cover of
$\pi_V(E)$. For each i, we may use the Besicovitch covering theorem
to find a constant $\xi$, depending only on n, and a family of balls
$\big(B_{ij}=B(x_{ij},r_{ij})\big)_{j\in\mathbb{N}}$ with
$r_{ij}=\frac {r_i} 2$ which is a $\delta$-cover of
$\pi_V^{-1}(B_i)\bigcap E$ such that
\begin{center}$\displaystyle\bigcup_jB(x_{ij},r_{ij})\subseteq
\pi_V^{-1}\Big(B(x_i,2r_i))\cap V\Big).$\end{center} Note that
$\widetilde{B}_i=B(x_i,2r_i)$. Then
\begin{eqnarray*}
\displaystyle\sum_i\mu_V\big(\widetilde{B}_i\big)^q\nu_V\big(\widetilde{B}_i\big)^t(4r_i)^s&\geq&
\xi^{-(q+t)}\displaystyle\sum_i(4r_i)^s\left(\sum_j\mu(B_{ij})\right)^q\left(\sum_j\nu(B_{ij})\right)^t \\
&\geq&\xi^{-(q+t)}\displaystyle\sum_{i,j}(4r_i)^s\mu(B_{ij})^q\nu(B_{ij})^t \\
&\geq&
4^s\xi^{-(q+t)}\displaystyle\sum_{i,j}\mu(B_{ij})^q\nu(B_{ij})^t(2r_{ij})^s.\end{eqnarray*}
Consequently, as $(B_i)_i$ way any $\delta$-cover of $\pi_V(E)$, we
conclude that
\begin{center}$\displaystyle\overline{\mathcal{H}}^{q,t,s}_{\mu,\nu,\delta}(E)\leq
4^{-s}\xi^{(q+t)}\displaystyle\overline{\mathcal{H}}^{q,t,s}_{\mu_V,\nu_V,2\delta}(\pi_V(E)).$\end{center}
Letting $\delta\downarrow0$, gives that
\begin{center}$\displaystyle\overline{\mathcal{H}}^{q,t,s}_{\mu,\nu}(E)\leq
4^{-s}\xi^{(q+t)}\displaystyle\overline{\mathcal{H}}^{q,t,s}_{\mu_V,\nu_V}(\pi_V(E)).$\end{center}
Thus
\begin{eqnarray*}
\overline{\mathcal{H}}_{\mu,\nu}^{q,t,s}(E) &\leq&4^{-s}\xi^{(q+t)}
\overline{\mathcal{H}}_{\mu_V,\nu_V}^{q,t,s}(\pi_V(E)) \\
&\leq&4^{-s}\xi^{(q+t)}\overline{\mathcal{H}}_{\mu_V,\nu_V}^{q,t,s}(\pi_V(\supp\mu))\\
&\leq&4^{-s}\xi^{(q+t)}\mathcal{H}_{\mu_V,\nu_V}^{q,t,s}(\supp\mu_V).
\end{eqnarray*}
The arbitrary on $E$ implies that
\begin{equation*}
\mathcal{H}_{\mu,\nu}^{q,t,s}(\supp\mu)\leq4^{-s}\xi^{(q+t)}
\mathcal{H}_{\mu_V,\nu_V}^{q,t,s}(\supp\mu_V)
\end{equation*}
and the result yields.
\begin{remark}
 Notice that in the case where $t=0$ or $q=0$, the preceding results were
treated by O$'$Neil in \cite{28}. Also, when $s=0$, Douzi and Selmi
investigated the projection properties of the mutual Hausdorff,
packing and pre-packing measures in \cite{4}. They derived global
bounds on the relative multifractal dimensions of a projection of a
measures in terms of its original relative multifractal dimensions.
\end{remark}
\section{Application}
This section is devoted to study the behavior of projections of
measures obeying to the mutual multifractal formalism. More
precisely, we prove that for
$(q,t)\in\Big\{\big(]-\infty,0[^2\big)\cup$\\$\big(]-\infty,0[\times
]0,1]\big)\cup\big(]0,1]\times ]-\infty,0[\big)\Big\}$, if the
mutual multifractal formalism holds for the couple $(\mu, \nu)$  at
$\alpha=-\frac{\partial B_{\mu,\nu}(q,t)}{\partial q}$ and
$\beta=-\frac{\partial B_{\mu,\nu}(q,t)}{\partial t}$, it holds for
$(\mu_V, \nu_V)$ for all $V\in G_{n,m}$. Before detailing our
results let us recall the mutual multifractal formalism introduced
by Svetova \cite{6, 8}. For $\alpha, \beta\geq 0$, let
 $$
E_{\mu,\nu}(\alpha,\beta)=\left\{x\in \supp\mu\cap \supp\nu \:;\:
\lim_{r\to 0}\frac{\log \big(\mu(B(x,r))\big)}{\log r}=\alpha
~\text{and}~\lim_{r\to 0}\frac{\log \big(\nu(B(x,r))\big)}{\log
r}=\beta \right\}.
 $$

We are interested to the estimation of the Hausdorff and packing
dimension of $E_{\mu,\nu}(\alpha,\beta)$. Let us mention that in the
last decay there has been a great interest for the multifractal
analysis and positive results have been written in various
situations (see for example \cite{2, 3, 10, 4, 25}). Our purpose in
the following theorem is to prove the result of Theorem 3 in
\cite{6} under less restrictive hypotheses.
\begin{theorem}\label{ah}
Let $\mu, \nu$ be two compactly supported Borel probability measures
on  $\mathbb{R}^n$. If $B_{\mu,\nu}$ is differentiable at $(q,t)$
and we set $\alpha=-\frac{\partial B_{\mu,\nu}(q,t)}{\partial q}$
and $\beta=-\frac{\partial B_{\mu,\nu}(q,t)}{\partial t}$. Assume
that
$\mathcal{H}_{\mu,\nu}^{q,t,B_{\mu,\nu}(q,t)}\big(\supp\mu\cap\supp\nu\big)>0$.
Then, we have
 $$
\dim_H E_{\mu,\nu}(\alpha,\beta)=\dim_P E_{\mu,\nu}(\alpha,\beta)=
B^*_{\mu,\nu}(\alpha,\beta)=b^*_{\mu,\nu}(\alpha,\beta),
 $$
where $f^*(\alpha,\beta)=\displaystyle\inf_{q,t}\big(\alpha q+ \beta
t +f(\alpha,\beta)\big)$ denotes the Legendre transform of the
function $f$. Here $\dim_H$ and $\dim_P$ denote the Hausdorff and
packing dimensions (see \cite{11} for the definitions) and we say
that the mutual multifractal formalism is valid.
\end{theorem}

\noindent{\bf Proof.} It is known (for instance, see \cite{6}) that,
for all  reals $\alpha$ and $\beta$, one has
$$
\dim_P E_{\mu,\nu}(\alpha,\beta)\leq \alpha q+\beta t+
B_{\mu,\nu}(q,t).
$$
\autoref{ah} is a consequence from the following lemmas.
\begin{lemma} \label{BBB}Let $\eta_1,\eta_2>0$. We set $\alpha=-\frac{\partial B_{\mu,\nu}(q,t)}{\partial q}$
and $\beta=-\frac{\partial B_{\mu,\nu}(q,t)}{\partial t}$. Then
\begin{eqnarray*}
{\mathcal H}^{\alpha q+\beta
t+B_{\mu,\nu}(q,t)-\eta_1-\eta_2}\big(E_{\mu,\nu}(\alpha,\beta)\big)\geq2^{\alpha
q+\beta t-\eta_1-\eta_2}{\mathcal
H}^{q,t,B_{\mu,\nu}(q,t)}_{\mu,\nu}\big(E_{\mu,\nu}(\alpha,\beta)\big).\end{eqnarray*}
\end{lemma}
\noindent{\bf Proof.} We treat the case $q\leq0$ and $t\le0$. The
other cases are proved similarly. The result is true for $q=t=0$, so
we may assume that $q < 0$ and $t<0$. \\For $m\in\mathbb{N}^\ast$,
write
\begin{multline*}
{E}_m:=\Big\{x\in E_{\mu,\nu}(\alpha,\beta)\:;\:\frac{\log
\big(\mu(B(x,r))\big)}{\log r}\le\alpha-\frac{\eta_1} q\\
~\text{and}~\frac{\log \big(\nu(B(x,r))\big)}{\log r}
\le\beta-\frac{\eta_2} t ~\text{for}~ 0<r<\frac 1 m\Big\}.
\end{multline*} Given $F\subseteq E_m$,
$0<\delta<\displaystyle\frac 1 m $ and $\Big(B(x_i,r_i)\Big)_i$ a
centered $\delta$-covering of $F$, we have
\begin{center}$\displaystyle\frac{\log\mu(B(x_i,r_i))}{\log
r_i}\leq\alpha-\frac{\eta_1} q $ and
$\displaystyle\frac{\log\nu(B(x_i,r_i))}{\log
r_i}\leq\beta-\frac{\eta_2} q .$ \end{center} Thus
\begin{center}$\mu(B(x_i,r_i))^q\leq {r_i}^{\alpha q-\eta_1}$
\text{and} $\nu(B(x_i,r_i))^t\leq {r_i}^{\beta
t-\eta_2}.$\end{center} Hence
\begin{center}$\mu(B(x_i,r_i))^q\nu(B(x_i,r_i))^t(2r_i)^{B_{\mu,\nu}(q,t)}\leq
2^{B_{\mu,\nu}(q,t)}{r_i}^{\alpha q+\beta
t+B_{\mu,\nu}(q,t)-\eta_1-\eta_2}.$ \end{center}So
\begin{eqnarray*}
   {\overline{\mathcal H}}_{\mu,\nu,\delta}^{q,t,B_{\mu,\nu}(q,t)}(F)&\le&
 \sum_i\mu(B(x_i,r_i))^q\nu(B(x_i,r_i))^t(2r_i)^{B_{\mu,\nu}(q)} \\
   &\le&2^{-\alpha q-\beta
t+\eta_1+\eta_2}\sum_i{(2r_i)}^{\alpha q+\beta
t+B_{\mu,\nu}(q,t)-\eta_1-\eta_2}.
\end{eqnarray*}
We can deduce that
$${\overline{\mathcal
H}}_{\mu,\nu,\delta}^{q,t,B_{\mu,\nu}(q,t)}(F)\leq2^{-\alpha q-\beta
t+\eta_1+\eta_2}\overline{{\mathcal H}}_{\delta}^{\alpha q+\beta
t+B_{\mu,\nu}(q,t)-\eta_1-\eta_2}(F).
$$
Letting $\delta\searrow0$ gives that

\begin{eqnarray*}
   {\overline{\mathcal
H}}_{\mu,\nu}^{q,t,B_{\mu,\nu}(q,t)}(F)&\le& 2^{-\alpha q-\beta
t+\eta_1+\eta_2}\overline{{\mathcal H}}^{\alpha q+\beta
t+B_{\mu,\nu}(q,t)-\eta_1-\eta_2}(F) \\
   &\le& 2^{-\alpha q-\beta
t+\eta_1+\eta_2}{\mathcal H}^{\alpha q+\beta
t+B_{\mu,\nu}(q,t)-\eta_1-\eta_2}(E_m)
\end{eqnarray*}
 for all $F\subseteq
E_m.$ Hence
$${\mathcal H}_{\mu,\nu}^{q,t,B_{\mu,\nu}(q,t)}(E_m)\leq2^{-\alpha q-\beta
t+\eta_1+\eta_2}{\mathcal H}^{\alpha q+\beta
t+B_{\mu,\nu}(q,t)-\eta_1-\eta_2}(E_m)$$ and the result follows
since $E_{\mu,\nu}(\alpha,\beta)=\displaystyle\bigcup_m E_m$.
\begin{lemma}\label{ll}
We have \;
$\mathcal{H}_{\mu,\nu}^{q,t,B_{\mu,\nu}(q,t)}\Big(\big(\supp\mu\cap
\supp\nu\big)\backslash E_{\mu,\nu}(\alpha,\beta)\Big)=0.$
\end{lemma}
\noindent{\bf Proof.} Let us introduce, for
$\alpha,\beta\in\mathbb{R}$
\begin{eqnarray*}\label{1}
F_{\alpha,\beta}=\left\{x;\;\limsup_{r\to 0}\frac{\log
\big(\mu(B(x,r))\big)}{\log r }>\alpha,\;\text{or}\; \limsup_{r\to
0}\frac{\log \big(\nu(B(x,r))\big)}{\log r}>\beta \right\},
\end{eqnarray*}
\begin{eqnarray*}\label{1}
F_{\alpha,\beta}^1=\left\{x;\;\liminf_{r\to 0}\frac{\log
\big(\mu(B(x,r))\big)}{\log r }<\alpha,\;\text{or}\; \liminf_{r\to
0}\frac{\log \big(\nu(B(x,r))\big)}{\log r}<\beta \right\},
\end{eqnarray*}
\begin{eqnarray*}\label{1}
F_{\alpha,\beta}^2=\left\{x;\;\limsup_{r\to 0}\frac{\log
\big(\mu(B(x,r))\big)}{\log r }>\alpha,\;\text{or}\; \liminf_{r\to
0}\frac{\log \big(\nu(B(x,r))\big)}{\log r}<\beta \right\},
\end{eqnarray*}
\begin{eqnarray*}\label{1}
F_{\alpha,\beta}^3=\left\{x;\;\liminf_{r\to 0}\frac{\log
\big(\mu(B(x,r))\big)}{\log r }<\alpha,\;\text{or}\; \limsup_{r\to
0}\frac{\log \big(\nu(B(x,r))\big)}{\log r}>\beta \right\}
\end{eqnarray*}
We have to prove that
\begin{eqnarray}\label{2}
\mathcal{H}_{\mu,\nu}^{q,t,B_{\mu,\nu}(q,t)}(F_{\alpha,\beta})=0\:\text{for
every}\: \alpha>-\frac{\partial B_{\mu,\nu}(q,t)}{\partial
q}\;\text{and}\; \beta>-\frac{\partial B_{\mu,\nu}(q,t)}{\partial t}
\end{eqnarray}
\begin{eqnarray}\label{3}
\mathcal{H}_{\mu,\nu}^{q,t,B_{\mu,\nu}(q,t)}(F_{\alpha,\beta}^1)=0\:\text{for
every}\: \alpha<-\frac{\partial B_{\mu,\nu}(q,t)}{\partial
q}\;\text{and}\; \beta<-\frac{\partial B_{\mu,\nu}(q,t)}{\partial t}
\end{eqnarray}
\begin{eqnarray}\label{4}
\mathcal{H}_{\mu,\nu}^{q,t,B_{\mu,\nu}(q,t)}(F_{\alpha,\beta}^2)=0\:\text{for
every}\: \alpha>-\frac{\partial B_{\mu,\nu}(q,t)}{\partial
q}\;\text{and}\; \beta<-\frac{\partial B_{\mu,\nu}(q,t)}{\partial t}
\end{eqnarray}
and
\begin{eqnarray}\label{5}
\mathcal{H}_{\mu,\nu}^{q,t,B_{\mu,\nu}(q,t)}(F_{\alpha,\beta}^3)=0\:\text{for
every}\: \alpha<-\frac{\partial B_{\mu,\nu}(q,t)}{\partial
q}\;\text{and}\; \beta>-\frac{\partial B_{\mu,\nu}(q,t)}{\partial t}
\end{eqnarray}
Let us sketch the proof of assertion \eqref{2}. Given
\begin{center}
$\alpha>-\frac{\partial B_{\mu,\nu}(q,t)}{\partial q}$ and
$\beta>-\frac{\partial B_{\mu,\nu}(q,t)}{\partial t}$, \end{center}
we can choose $h>0$ such that
\begin{center}
$B_{\mu,\nu}(q-h,t)<B_{\mu,\nu}(q,t)+\alpha h$ and
$B_{\mu,\nu}(q,t-h)<B_{\mu,\nu}(q,t)+\beta h$.
\end{center}
Which implies,
$$
{\mathcal P}^{q-h,t,B_{\mu,\nu}(q,t)+\alpha
h}_{\mu,\nu}\big(\supp\mu\cap \supp\nu\big)=0
$$
and
$$
{\mathcal P}^{q,t-h,B_{\mu,\nu}(q,t)+\beta
h}_{\mu,\nu}\big(\supp\mu\cap \supp\nu\big)=0.
$$
Let $\delta>0$. For each $x\in F_{\alpha,\beta}$, there exists
$0<r_x<\delta$ such that
 $$
\mu(B(x,r_x))\leq r_x^\alpha\qquad\text{or}\qquad\nu(B(x,r_x))\leq
r_x^\beta.
 $$
The family $\Big(B(x,r_x)\Big)_{x\in F_{\alpha,\beta}}$ is then a
centered $\delta$-covering of $ F_{\alpha,\beta}.$ Using
Besicovitch's Covering Theorem, we can construct $\xi$ finite or
countable sub-families
$\Big(B(x_{1j},r_{1j})\Big)_j,...,\Big(B(x_{\xi j},r_{\xi
j})\Big)_j$ such that each  $
F_{\alpha,\beta}\subseteq\displaystyle\bigcup_{i=1}^\xi\bigcup_jB(x_{ij},r_{ij})$
and $\Big(B(x_{ij},r_{ij})\Big)_j$ is a $\delta$-packing of
$F_{\alpha,\beta}.$ Observing that
 $$
\mu(B(x_{ij},r_{ij}))^q\nu(B(x_{ij},r_{ij}))^t(2r_{ij})^{B_{\mu,\nu}(q,t)}\leq
\mu(B(x_{ij},r_{ij}))^{q-h} \nu(B(x_{ij},r_{ij}))^t
(2r_{ij})^{B_{\mu,\nu}(q,t)+\alpha h}
 $$
or
 $$
\mu(B(x_{ij},r_{ij}))^q\nu(B(x_{ij},r_{ij}))^t(2r_{ij})^{B_{\mu,\nu}(q,t)}\leq
\mu(B(x_{ij},r_{ij}))^{q} \nu(B(x_{ij},r_{ij}))^{t-h}
(2r_{ij})^{B_{\mu,\nu}(q,t)+\beta h},
 $$
we obtain
 $$
\overline{\mathcal{H}}_{\mu,\nu}^{q,t,B_{\mu,\nu}(q,t)}(F_{\alpha,\beta})
\leq\xi\overline{\mathcal{P}}_{\mu,\nu}^{q-h,t,B_{\mu,\nu}(q,t)+\alpha
h}( F_{\alpha,\beta})
 $$
or
$$
\overline{\mathcal{H}}_{\mu,\nu}^{q,t,B_{\mu,\nu}(q,t)}(F_{\alpha,\beta})
\leq\xi\overline{\mathcal{P}}_{\mu,\nu}^{q,t-h,B_{\mu,\nu}(q,t)+\beta
h}( F_{\alpha,\beta}).
 $$
Remark that, in the last inequality, we can replace
$F_{\alpha,\beta}$ by any arbitrary subset of $F_{\alpha,\beta}.$
Then, we can finally conclude that
\begin{eqnarray*}
\mathcal{H}_{\mu,\nu}^{q,t,B_{\mu,\nu}(q,t)}(
F_{\alpha,\beta})&\leq&\xi\mathcal{P}_{\mu,\nu}^{q-h,t,B_{\mu,\nu}(q,t)+\alpha
h}( F_{\alpha,\beta})\\&\leq& \xi {\mathcal
P}^{q-h,t,B_{\mu,\nu}(q,t)+\alpha h}_{\mu,\nu}\big(\supp\mu\cap
\supp\nu\big)=0
\end{eqnarray*}
 or
\begin{eqnarray*}
\mathcal{H}_{\mu,\nu}^{q,t,B_{\mu,\nu}(q,t)}(
F_{\alpha,\beta})&\leq&\xi\mathcal{P}_{\mu,\nu}^{q,t-h,B_{\mu,\nu}(q,t)+\beta
h}( F_{\alpha,\beta})\\&\leq& \xi {\mathcal
P}^{q,t-h,B_{\mu,\nu}(q,t)+\beta h}_{\mu,\nu}\big(\supp\mu\cap
\supp\nu\big)=0.
\end{eqnarray*} The proof of \eqref{3}, \eqref{4}
and \eqref{5} is similar to \eqref{2}.$\hfill\square$

 \bigskip
Let us return to the proof of \autoref{ah}. By \autoref{BBB} and
\autoref{ll}, we have for all $\eta_1, \eta_2>0$,
$${\mathcal H}^{\alpha
q+\beta
t+B_{\mu,\nu}(q,t)-\eta_1-\eta_2}\big(E_{\mu,\nu}(\alpha,\beta)\big)\geq2^{\alpha
q+\beta t-\eta_1-\eta_2}{\mathcal
H}^{q,t,B_{\mu,\nu}(q,t)}_{\mu,\nu}\big(E_{\mu,\nu}(\alpha,\beta)\big)>0.$$
Whence,
$$
\dim_H E_{\mu,\nu}(\alpha,\beta)\geq \alpha q+\beta
t+B_{\mu,\nu}(q,t)-\eta_1-\eta_2.
$$
Letting $\eta_1\to 0$ and $\eta_2\to 0$ yields
$$
\dim_H E_{\mu,\nu}(\alpha,\beta)\geq \alpha q+\beta
t+B_{\mu,\nu}(q,t).
$$
Which achieves  the proof of \autoref{ah}.$\hfill\square$

 \bigskip
In the following we study the validity of the multifractal formalism
under projection.
\begin{theorem}\label{x} Let $\mu, \nu$ be two compactly supported Borel
probability measures on $\mathbb{R}^n$ such that
$\supp\mu=\supp\nu$. For
$(q,t)\in\big(]-\infty,0[^2\big)\cup\big(]-\infty,0[\times
]0,1]\big)\cup\big(]0,1]\times ]-\infty,0[\big)$, suppose that

\begin{hyp}[$H_1$]
$\mathcal{H}_{\mu,\nu}^{q,t,B_{\mu,\nu}(q,t)}(\supp\mu\cap\supp\nu)>0,$
\end{hyp}

\begin{hyp}[$H_2$]
$B_{\mu,\nu}$ is differentiable at $(q,t)$,
\end{hyp}
\noindent Then, for all $ V \in G_{n,m}$, we have
\begin{eqnarray*}
  \dim_H E_{\mu_V,\nu_V}\big(\alpha,\beta\big)&=&
  \dim_P E_{\mu_V,\nu_V}\big(\alpha,\beta\big)
  =\dim_H E_{\mu,\nu}\big(\alpha,\beta\big)\\&=&
  \dim_P E_{\mu,\nu}\big(\alpha,\beta\big)
  =
  B^*_{\mu,\nu}\big(\alpha,\beta\big)=b^*_{\mu,\nu}\big(\alpha,\beta\big).
\end{eqnarray*}
where $\alpha=-\frac{\partial B_{\mu,\nu}(q,t)}{\partial q}$ and
$\beta=-\frac{\partial B_{\mu,\nu}(q,t)}{\partial t}$.
\end{theorem}
\noindent{\bf Proof.} By using \autoref{j}, \autoref{k} and $(H_1)$,
we have
\begin{equation}\label{qq}
b_{\mu,\nu}(q,t)=B_{\mu,\nu}(q,t)=b_{\mu_V,\nu_V}(q,t)=B_{\mu_V,\nu_V}(q,t),
\qquad\forall V\in G_{n,m}.
\end{equation}
$(H_1)$, \eqref{www} and \eqref{qq} ensure that
 $$
\mathcal{H}_{\mu_V,\nu_V}^{q,t,B_{\mu_V,\nu_V}(q,t)}(\supp\mu_V)\ge
\mathcal{H}_{\mu,\nu}^{q,t,B_{\mu,\nu}(q,t)}(\supp\mu)>0,
\qquad\forall V\in G_{n,m}.
 $$
So, \autoref{ah} and the equalities \eqref{qq} imply that
\begin{eqnarray*}
\dim_H {E}_{\mu_V,\nu_V}\big(\alpha,\beta\big)\geq \alpha q+\beta t
+B_{\mu,\nu}(q,t).
\end{eqnarray*}
The other estimation is satisfied since
\begin{eqnarray*}
\dim_P{E}_{\mu_V,\nu_V}\big(\alpha,\beta\big)&\leq& \alpha q+ \beta
t +B_{\mu_V,\nu_V}(q,t)\\& = &\alpha q+\beta t +B_{\mu,\nu}(q,t).
\end{eqnarray*}
Which achieves the proof of \autoref{x}.$\hfill\square$
\begin{remark}
Let $\mu$ and $\nu$ be two compactly supported Borel probability
measures on $\mathbb{R}^n$. We write for $\gamma\geq0$,
 $$
\mathcal{B}_{\mu,\nu}(\gamma)=\left\{x\in \supp\mu\cap \supp\nu ;\:
\:\lim_{r\to 0}\frac{\log \big(\mu(B(x,r))\big)}{\log\big(
\nu(B(x,r)\big)}=\gamma\right\}.
 $$
It is clear that
 $$
\displaystyle
\bigcup_{\displaystyle{\substack{(\alpha,\beta)\in\mathbb{R}_+\times
\mathbb{R}_+^*, \\\frac{\alpha}{\beta}=\gamma}}}\;
E_{\mu,\nu}\big(\alpha,\beta\big) \subset
\mathcal{B}_{\mu,\nu}(\gamma).
 $$
The latter union is composed by an uncountable number of pairwise
disjoint nonempty sets. \autoref{ah} shows that surprisingly the
Hausdorff and packing dimension of $\mathcal{B}_{\mu,\nu}(\gamma)$
is fully carried by some subset $E_{\mu,\nu}\big(\alpha,\beta\big)$.
Together with \autoref{x}, this relationship provides a lower bound
to the relative multifractal spectra of the projections of a
measures introduced in Theorem 4.2 in \cite{4}.
\end{remark}


\begin{thebibliography}{99}
\bibitem{NBC} N. Attia, B. Selmi and Ch.
Souissi. {\it Some density results of relative multifractal
analysis}. Chaos, Solitons and Fractals. (2017), vol.~103, pp.~1-11.

\bibitem{BSB} L. Barreira, B. Saussol and J. Schmeling. {\it Higher-dimensional multifractal
analysis}. J. Math. Pures Appl., (2002), Vol.~81, pp.~ 67-91.

\bibitem{BP} L. Barreira and P. Doutor. {\it Birkhoff Averages for Hyperbolic Flows:
Variational Principles and Applications}. Journal of Statistical
Physics. (2004), Vol.~115, pp.~ 1567-1603.

\bibitem{BP1} L. Barreira and P. Doutor. {\it Almost additive multifractal analysis}.
J. Math. Pures Appl., (2009), Vol.~92, pp.~ 1-17.

\bibitem{BP2} L. Barreira and P. Doutor. {\it Dimension spectra of almost additive sequences}.
Nonlinearity. (2009), Vol.~22, pp.~ 2761-2773.

\bibitem{BCW} L. Barreira, Y. Cao and J. Wang. {\it Multifractal Analysis of Asymptotically Additive
Sequences}. J. Stat. Phys., (2013), Vol.~153, pp.~ 888-910.

\bibitem{1} J.Barral and I.Bhouri. {\it Multifractal analysis for projections
of Gibbs and related measures}. Ergodic Theory and Dynamic systems.
(2011), vol.~31, pp.~673-701.

\bibitem{BF} J. Barral and D.J. Feng. \emph{Projections of planar Mandelbrot
measures}. \href{http://arXiv.org/archive/math\#arXiv 1501protecd
kern+.2777em elax00875v1.}{ arXiv : 1605.09083v1}, (2016).

\bibitem{2} F. Ben Nasr. {\it Analyse multifractale de mesures}. CR Acad Sci Paris
Ser I. (1994), vol.~319, pp.~807-10.

\bibitem{3}F. Ben Nasr, I. Bhouri and Y. Heurteaux.{\it The validity of the
multifractal formalism: results and examples}. Adv. in Math. (2002),
vol.~165, pp.~264-284.

\bibitem{10} J. Cole. {\it Relative multifractal analysis}. Chaos, Solitons and Fractals.
(2000), vol.~11, pp.~2233-2250.

\bibitem{D1} C. Dai, Y. Li. {\it A multifractal formalism in a probability
space}. Chaos Solitons Fractals. 2006, vol.~27, pp.~57-73.

\bibitem{D3} C. Dai, Y. Li. {\it Multifractal dimension inequalities in a
probability space}. Chaos Solitons Fractals. (2007), vol.~34,
pp.~213-223.

\bibitem{D4} M. Dai, X. Peng. and W. Li. {\it Relative Multifractal Analysis in a
Probability Space}. Int. J. Nonlinear Sci., (2010), vol.~10,
pp.~313-319.

\bibitem{D} M. Dai. \emph{Mixed self-conformal multifractal measures}.
Analysis in Theory and Applications. {\bf25} (2009), 154-165.

\bibitem{D11} M. Dai and Y. shi. \emph{Typical behavior of mixed $L^q$-dimensions}.
Nonlinear Analysis: Theory, Methods $\&$ Applications. {\bf72}
(2010), 2318-2325.

\bibitem{D21} M. Dai and W. Li. \emph{The mixed $L^q$-spectra of self-conformal measures satisfying the weak
separation condition}. J. Math. Anal. Appl., {\bf382} (2011),
140-147.

\bibitem{D31} M. Dai, C. Wang and H. Sun. \emph{Mixed generalized dimensions of random self-similar
measures}. Int. J. Nonlinear. Sci., {\bf13} (2012), 123-128.

\bibitem{D41} M. Dai, J. Houa, J. Gaob, W. Suc, L. Xid and D. Ye. \emph{Mixed multifractal analysis of China and US stock index
series}.  Chaos, Solitons $\&$ Fractals. {\bf87} (2016), 286-275.

\bibitem{D51} M. Dai, S. Shao, J. Gao, Y. Sun and W. Su. \emph{Mixed multifractal analysis of crude oil, gold and exchange rate series}.
 Fractals. {\bf24} (2016), 1-7.

\bibitem{4} Z. Douzi and B. Selmi.{\it Multifractal variation for projections
of measures}. Chaos, Solitons and Fractals. (2016), vol.~91,
pp.~414-420.

\bibitem{11} K. J. Falconer. {\it The Geometry of Fractal sets}. Cambridge
univ. Press New. York-London. (1985), vol.~85.

\bibitem{12} K. J. Falconer, J. D. Howroyd. {\it Packing Dimensions of
Projections and Dimensions Profiles}. Math. Proc. Cambridge Philos.
Soc. (1997), vol.~121, pp.~269-286.

\bibitem{13} K. J. Falconer, P. Mattila. {\it The Packing Dimensions of
Projections and Sections of Measures}. Math. Proc. Cambridge Philos.
Soc. (1996), vol.~119, pp.~695-713.


\bibitem{15} K. J. Falconer, X. Jim. {\it Exact dimensionality and
projections of random self-similar measures and sets}. J.Lond.Math.
Soc. (2014), vol.~90, pp.~388-412.

\bibitem{16} M. Hochman, P. Shmerkin. {\it Local entropy averages and projections of
fractal measures}. \href{http://arXiv.org/archive/math\#arXiv;
0910protecd kern+.2777em elax1956v1}{arXiv : 0910.1956v1}, (2009).


\bibitem{19} R. Kaufman. {\it On Hausdorff dimension of projections}.
Mathematika. (1968), vol.~15, pp.~153-155.

\bibitem{21} J .M. Marstrand. Some fundamental geometrical properties
 of plane sets of fractional dimensions. Proceedings of the London
Mathematical Society. (1954), vol.~4, pp.~257-302.

\bibitem{22} P. Mattila. {\it Hausdorff dimension, orthogonal projections
and intersections with planes.} Annales Academiae Scientiarum
Fennicae. Series A I. Mathematica. (1975), vol.~1, pp.~227-244.

\bibitem{23} P. Mattila. {\it The Geometry of Sets and Measures in Euclidean
 Spaces}. Cambridge University Press, Cambrdige. (1995).

\bibitem{LOW} J. Li, L. Olsen and M. Wu. \emph{Bounds for the $L^q$-spectra of self-similar measures without any
separation conditions}. J. Math. Anal. Appl., {\bf 387} (2012),
77-89.

\bibitem{25} L. Olsen. {\it A multifractal formalism}. Advances in Mathematics.
(1995), vol.~116, pp.~82-196.

\bibitem{Ol1} L. Olsen. \emph{Mixed generalized dimensions of self-similar measures}.
J. Math. Anal. Appl., {\bf 306} (2005), 516-539.

\bibitem{Ol2} L. Olsen. \emph{Bounds for the $L^q$-spectra of a self-similar multifractal not satisfying
the open set condition}. J. Math. Anal. Appl., {\bf 355} (2005),
12-21.

\bibitem{Ol} L. Olsen. \emph{On the inverse multifractal formalism}. Glasgow Mathematical Journal. {\bf 52} (2010),
179-194.

\bibitem{27} T. C. O'Neil. {\it The multifractal
spectrum of quasi self-similar measures}. Journal of Mathematical
Analysis and its Applications. (1997), vol.~211, pp.~233-257.

\bibitem{28} T. C. O'Neil. {\it The multifractal spectra of projected measures in Euclidean
spaces}. Chaos, Solitons and Fractals. (2000), vol.~11, pp.~901-921.

\bibitem{JP} J. Peyri\'{e}re. {\it Multifractal measures}. In: Proceedings of the NATO Advanced Study
Institute on Probabilistic and Stochastic Methods in Analysis with
Applications. Il Ciocco, NATO ASI Series, Series C: Mathematical and
physical sciences, vol. 372, Kluwer Academic Press, Dordrecht,
(1992), pp. 175-186.

\bibitem{B} B. Selmi. {\it Multifractal dimensions for projections of
 measures}. Preprint. 2017.

\bibitem{SB} B. Selmi. \emph{A note on the effect of projections on both measures
 and the generalization of $q$-dimension capacity}. Probl. Anal. Issues Anal., (2016), vol.~5, pp.~38 - 51.

\bibitem{SBS} B. Selmi and N.Yu. Svetova. \emph{On the projections of mutual $L^{q,t}$-spectrum}. Probl. Anal. Issues Anal., (2017), vol.~6,
pp.~94 - 108.

\bibitem{31} P. Shmerkin. {\it Projections of self-similar and related fractals: Asurvey of recent
developments.}  \href{http://arXiv.org/archive/math\#arXiv
1501protecd kern+.2777em elax00875v1.}{ arXiv : 1501.00875v1},
(2015).

\bibitem{32} P.Shmerkin, B. Solomyak. {\it Absolute continuity of self-similar measures, their projections and
convolution,}  \href{http://arXiv.org/archive/math\#arXiv 1406
protecd kern+.2777em elax0204v1.}{arXiv : 1406.0204v1}, (2014).

 \bibitem{SSSS}  M. Slimane. \emph{Baire typical results for mixed H\"{o}lder spectra on product of
 continuous besov or oscillation spaces}. Mediterr. J. Math., {\bf 13}
(2016), 1513-1533.

\bibitem{5} N.Yu. Svetova. {\it Conditional and mutual multifractal spectra.
Definition and basic properties}. Tr. Petrozavodsk. Gos. Univ. Ser.
Mat., (2003), vol.~10, pp.~41-58.

\bibitem{6} N.Yu. Svetova. {\it Mutual multifractal spectra I: Exact
spectra}. Tr. Petrozavodsk. Gos. Univ. Ser. Mat., (2004), vol.~11,
pp.~41-46.

\bibitem{7} N.Yu. Svetova. {\it Mutual multifractal spectra II: Legendre and
Hentschel-Procaccia spectra, and spectra defined for partitions}.
Tr. Petrozavodsk. Gos. Univ. Ser. Mat., (2004), vol.~11, pp.~47-56.

\bibitem{8} N.Yu. Svetova. {\it An estimate for exact mutual multifractal
spectra}. Tr. Petrozavodsk. Gos. Univ. Ser. Mat., (2008), vol.~14,
pp.~59-66.

\bibitem{S} N.Yu. Svetova. {\it The property of convexity of mutual
multifractal dimension}. Tr. Petrozavodsk. Gos. Univ. Ser. Mat.,
(2010), vol.~17, pp.~15-24.

\end{thebibliography}
\end{document}